\documentclass[12pt]{article}

\usepackage{amsmath,amsfonts,amssymb,amsbsy,graphicx,xcolor}

\usepackage{color}

 \newcounter{definition}

\newcounter{sect1}
\newcounter{subsect1}

\newfont{\gotikai}{eufm10}

\newcommand {\beq}{\begin{equation}}
\newcommand {\eeq}{\end{equation}}

\newtheorem{Theorem}{Theorem}
\newtheorem{Lemma}{Lemma}
\newtheorem{Prop}{Proposition}

\newtheorem{Remark}{Remark}

\newtheorem{\theRemark}{\arabic{Remark}}
\newtheorem{definition}{\indent Definition}

\def\({\left(}
\def\){\right)}

\def\({\left(}
\def\){\right)}

\def\[{\left[}
\def\]{\right]}

\def\|{\left|}
\def\|{\right|}

\def\0{{\boldsymbol{0}}}

\def\DD{\textbf{Proof}}

\def\C{\mathbb C}

\def\DD{\textbf{Proof}}

\def\C{\mathbb C}

\def\00{{\boldsymbol{0}}}

\def\DD{\textbf{Proof}}

\def\C{\mathbb C}

\begin{document}

\begin{center}
\textbf{Limit theorems for the Wiener process with resetting}
\end{center}

\ \ \

\begin{center}
\textbf{A.V. Logachov, O.M.Logachova, A.A.Yambartsev, and K.A. Zaykov}
\end{center}

\ \ \

\textbf{Abstract}. We establish a large deviation principle for the trajectories of Wiener processes subject to random resets to the origin occurring according to a Poisson process. In addition to the pathwise large deviation principle, we identify the rate function and establish a large deviation principle for the supremum of the process over long time intervals.

\ \ \

\section{Introduction}

In this paper, we continue our investigation into large deviations for Markov processes with catastrophes \cite{LY0,LY3,LY2}. Within this framework, stochastic processes with resettings can be viewed as a very special -- yet remarkably prominent -- subclass of processes with catastrophes. Although technically a special case, reset-type dynamics stand out due to their rich structure and broad applicability in search processes, queueing systems, biology, and even in quantum physics. The volume of literature on resetting is substantial, and a comprehensive review is beyond our scope. We restrict ourselves to mentioning a few representative and survey-style works. 

The importance of resetting in modern stochastic modeling is reflected in the rapidly growing literature and, in particular, in the special issue “Stochastic Resetting: Theory and Applications” (Guest Editors: Anupam Kundu and Shlomi Reuveni; see their preface \cite{KR}), dedicated to the 10th anniversary of the foundational paper “Diffusion with Stochastic Resetting” \cite{EM}. Among the contributions to this issue, we refer to \cite{NG} for a review of resetting-type processes in particle systems, and to \cite{Z} for a comprehensive overview of rigorous mathematical results on resetting. Before this special issue, a valuable review of stochastic processes with resetting is provided in \cite{EMS}; see also the brief review \cite{GJ} for a concise introduction to the topic. For readers interested in reset-type processes in queueing systems, we refer to \cite{RBP}.

Large deviations for stochastic processes with resetting remain a relatively underexplored area of research. Notable exceptions include the paper \cite{MST} (see also \cite{Z}, the letter \cite{HT}, and the recent lecture notes \cite{ES}), where large deviation principles are established for integral (additive) functionals of diffusion processes with resetting. In \cite{CO}, large deviation techniques are employed to design stochastic resetting protocols that increase the probability of reaching desired configurations, such as promoting homogeneity in systems with network heterogeneity. It is also worth mentioning the paper \cite{M}, where large deviations were studied for the empirical measure, flow, and trajectories between successive resettings.

We are interested in large deviations at the level of process trajectories. In \cite{LY0}, we studied the Wiener process with catastrophes, where, at each Poissonian resetting time, a new position is chosen randomly according to a conditional distribution that depends on the state of the process just before the reset. We focused on a class of such distributions referred to as \textit{quasi-uniform catastrophes}, where the new state is selected (quasi-)uniformly between the previous state and zero. While our original goal was to establish a full large deviation principle, the results of that work were limited to a local large deviation principle in the functional space of cádlág functions, from which we were nevertheless able to derive the explicit form of the rate functional.

Here, we prove the (full) large deviation principle for the Wiener processes with resetting at the origin. In the literature on processes with catastrophes, the resetting at the origin corresponds to the so-called total catastrophes. We prove the principle in the metric space $\mathbb{L}[0,1]$, the space of all real-valued functions defined on the interval $[0,1]$ that are integrable in the sense of Lebesgue, with distance given by the integral of the absolute difference. Moreover, we prove the principle also for the extreme statistics.  

It appears that for stochastic processes with catastrophes, one could also hope to establish a large deviation principle in the metric space $\mathbb{L}[0,1]$. However, obtaining such results is considerably more difficult. For example, in \cite{LY0} the lower bound in the local large deviation principle becomes finite for certain continuous functions that are not absolutely continuous. This already indicates that the form of the rate functional in the large deviations setting for the Wiener process with catastrophes is significantly more involved, and that the proof itself becomes much more delicate.

Finally, it is worth noting that the classical approach developed by Feng and Kurtz \cite{Feng} for proving large deviation principles in Markov processes relies on certain uniform integrability conditions on the jump sizes, often referred to as the Cramér condition. Intuitively, this condition requires that the process does not make very large jumps, regardless of its current state. In the case of a Wiener process with resettings to the origin, the size of each jump depends on the current position of the process, which can become arbitrarily large over time. Although the expected jump size is finite at each point, it is not uniformly bounded across all states. As a result, the Cramér condition fails, and the standard large deviation framework of Feng and Kurtz does not directly apply. This motivates a more tailored, hands-on approach to establishing the large deviation principle in this setting.

The paper is organized as follows. In the next section, Section~\ref{dmr}, we recall some definitions and state the main results, Theorem~\ref{t.1} and Theorem~\ref{t.2}. The proofs of these theorems are given in Section~\ref{pmr}. Finally, in Section~\ref{Aux}, we collect some auxiliary lemmas. 

\section{Definitions and main results}\label{dmr}

Let $(\Omega,\mathfrak{F},\mathbf{P})$ be a probability space supporting a Wiener process $w(t)$ and a Poisson process $\nu(t)$ independent of it, with rate $\lambda$,  $\mathbf{E}\nu(1)=\lambda>0$.

\begin{figure}
    \centering
    \includegraphics[width=1\linewidth]{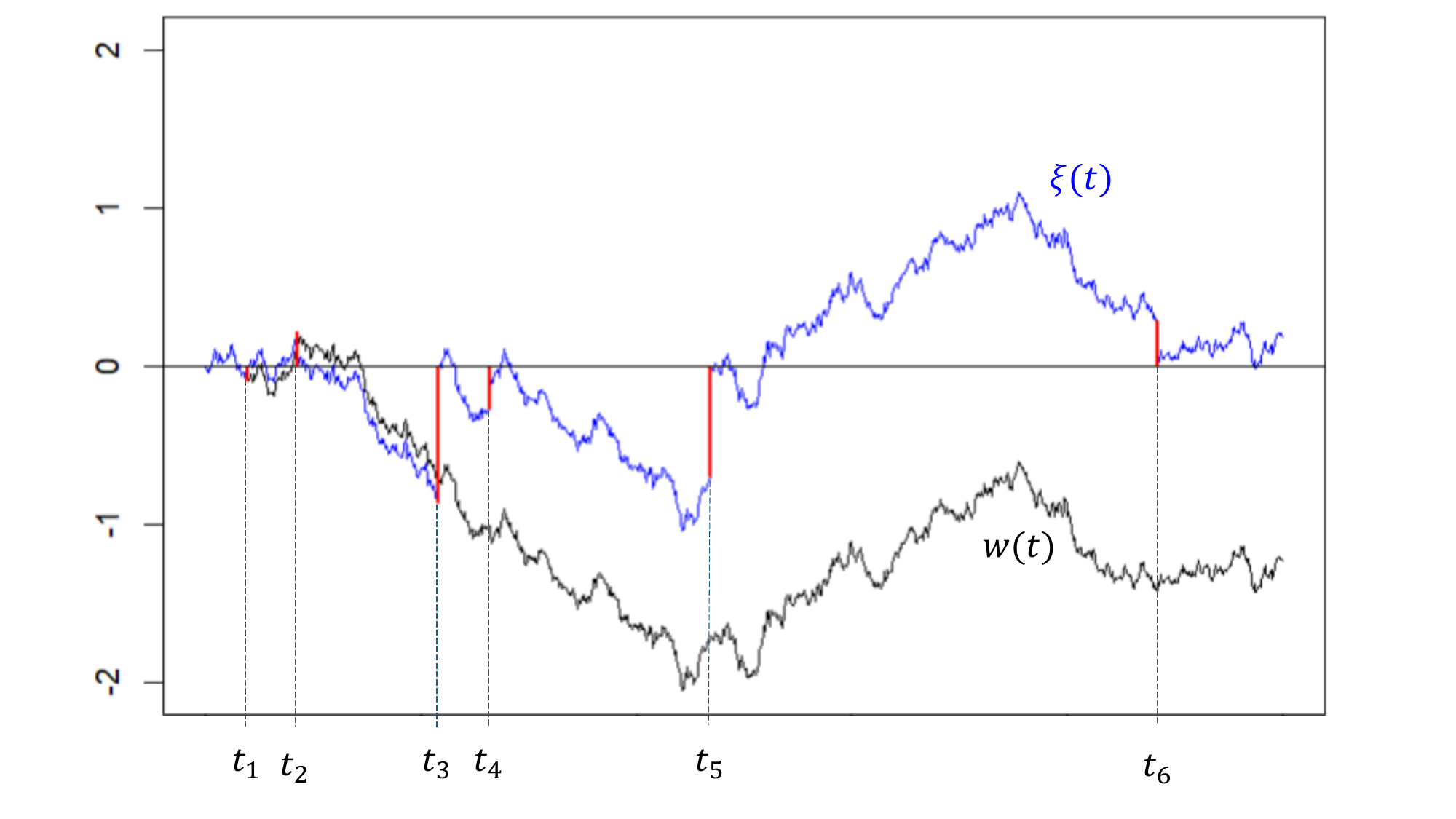}
    \caption{The black curve represents the original trajectory of the Wiener process $w$ over the time interval $[0,1]$; the blue curve shows the trajectory of the process $\xi$, with red vertical segments indicating the resetting jumps.}
    \label{fig1}
\end{figure}

The Wiener process with resetting is constructed as a stochastic process defined recursively as follows:
$$
\xi(t) := w(t) - \sum\limits_{k=1}^{\nu(t)} \xi(t_k-),
$$
where $t_1, t_2, \dots, t_k, \dots$ are the jump times of the Poisson process $\nu$; here and throughout, we assume that $\sum_{k=1}^0 = 0$ and $t_0 := 0$. Observe that the stochastic process $\xi$ can be written as
\begin{equation}\label{16-05-25-1-1}    
\xi(t)=w(t)-w(t_{\nu(t)}), \ \ \ t\in[0,1].
\end{equation}
See Figure~\ref{fig1} for the illustration of two representations of the process $\xi$.

We are interested in the large deviation behavior of the sequence
$$
\xi_n(t) := \frac{\xi(nt)}{n} = w_n(t) - \frac{1}{n} \sum\limits_{k=1}^{\nu(nt)} \xi(t_k-), \quad t \in [0,1],
$$
where $w_n(t) := \frac{w(nt)}{n}$. Representation \eqref{16-05-25-1-1} yields the following expression for the stochastic process $\xi_n$:
\begin{equation}\label{16-05-25-1-2}    
\xi_n(t)=w_n(t)-\frac{1}{n}w(t_{\nu(nt)})=w_n(t)-w_n\left(\frac{t_{\nu(nt)}}{n}\right), \ \ \ t\in[0,1].
\end{equation}

We consider the trajectories of the sequence of stochastic processes $\xi_n$ in the metric space $(\mathbb{L}[0,1], \rho_L)$, where
$$
\mathbb{L}[0,1] := \left\{ f : \int_0^1 |f(t)| \, dt < \infty \right\}, \quad
\rho_L(f,g) := \int_0^1 |f(t) - g(t)| \, dt, \quad f,g \in \mathbb{L}[0,1].
$$
We will also sometimes use the uniform metric as an auxiliary metric:
$$
\rho(f,g) := \sup\limits_{t \in [0,1]} |f(t) - g(t)|.
$$
For a set $A \subseteq \mathbb{L}[0,1]$, we denote
$$
\rho_L(f,A) := \inf\limits_{g \in A} \rho_L(f,g), \qquad \rho(f,A) := \inf\limits_{g \in A} \rho(f,g).
$$

Recall that in the metric space $(\mathbb{L}[0,1], \rho_L)$, functions that differ only on a set of Lebesgue measure zero (i.e., elements of the same equivalence class) are identified: the distance between any such functions is zero. Accordingly, we shall always work with a typical representative of each equivalence class, without specifying it explicitly.

For instance, when we say that a function $f \in \mathbb{L}[0,1]$ is continuous, we mean that there exists a continuous function $\hat{f}$ such that $\rho_L(f, \hat{f}) = 0$. Similarly, the statement that a function $f$ has bounded total variation on $[0,1]$ is understood as the existence of a function $\hat{f}$ with bounded total variation such that $\rho_L(f, \hat{f}) = 0$, and analogous conventions will be used throughout.

\medskip

For the reader’s convenience, we recall the necessary definitions.

\begin{definition} \label{d1.1}
A family of random elements $\xi_n$ satisfies the local large deviations principle (LLDP) in the metric space $(\mathbb{X}, \rho_X)$ with a rate function 
$I = I(x)\,:\, \mathbb{X} \rightarrow [0,\infty]$ and a normalizing function
$\psi(n)$ such that $\lim\limits_{n \rightarrow \infty} \psi(n) = \infty$,
if the following equality holds for any $x \in \mathbb{X}$:
$$
\lim\limits_{\varepsilon \rightarrow 0} \limsup\limits_{n \rightarrow \infty} \frac{1}{\psi(n)}
\ln \mathbf{P}(\rho_X(\xi_n, x) < \varepsilon)
=
\lim\limits_{\varepsilon \rightarrow 0} \liminf\limits_{n \rightarrow \infty} \frac{1}{\psi(n)}
\ln \mathbf{P}(\rho_X(\xi_n, x) < \varepsilon) = -I(x).
$$
\end{definition}

\begin{definition} \label{d1.2}
A family of random variables $\xi_n$ is said to be exponentially tight (ET) in the metric space $(\mathbb{X}, \rho_X)$ if, for any $N > 0$, there exists a compact set $K_N \subset \mathbb{X}$ such that
$$
\limsup_{n \rightarrow \infty} \frac{1}{\psi(n)} \ln \mathbf{P}(\xi_n \in (K_N)^c ) \leq -N,
$$
where $(K_N)^c$ denotes the complement of the set $K_N$.
\end{definition}

We denote the closure and interior of the set $B$ by $[B]$ and $(B)$, respectively.

\begin{definition} \label{d1.3}
  A family of random variables $\xi_n$ satisfies the LDP on the metric space $(\mathbb{X},\rho_X)$ with a \emph{good} rate function
$I = I(f)\,:\, \mathbb{X} \rightarrow [0,\infty]$ and normalizing function
$\psi(n)$, with $\lim\limits_{n\rightarrow\infty}\psi(n) = \infty$,
if for any $c \geq 0$, the set $\{ x \in \mathbb{X}\,:\, I(x) \leq c \}$ is compact,
and for every Borel set $B \in \mathfrak{B}(\mathbb{X})$ the following inequalities hold:
$$
 \limsup_{n \rightarrow \infty} \frac{1}{\psi(n)} \ln \mathbf{P}(\xi_n \in B ) \leq - I([B]),
$$
$$
\liminf_{n \rightarrow \infty} \frac{1}{\psi(n)} \ln \mathbf{P}(\xi_n \in B )\geq -I((B)),
$$
where $\mathfrak{B}(\mathbb{X})$ is the Borel $\sigma$-algebra on $\mathbb{X}$, $I(B) = \inf\limits_{x \in B} I(x)$,
and $I(\emptyset) = \infty$.
\end{definition}

In what follows, we will also use the following notation:
\begin{itemize}
\item $\mathbb{D}[0,1]$ is the space of c\`adl\`ag functions (i.e., right-continuous with left limits on the interval $[0,1]$);
\item $\mathbb{C}[0,1]$ is the space of continuous functions on $[0,1]$;
\item $\mathbb{AC}_0[0,1]$ is the space of absolutely continuous functions on $[0,1]$ satisfying $f(0)=0$;
\item $\widetilde{\mathbb{AC}}_0[0,1]$ is the set of functions $f\in \mathbb{D}[0,1]$ satisfying:
\begin{enumerate}
  \item $f(0)=0$;
  \item If $t^*$ is a point of discontinuity of $f$, then $f(t^*)=0$;
  \item If $t^*<t^{**}$ are two consecutive points of discontinuity (i.e., there are no other discontinuities between them), then $f$ is absolutely continuous on the interval $[t^*,t^{**})$;
\end{enumerate}
\item $\mathbb{V}[0,1]$ is the space functions with a finite variation. 
\end{itemize}

For a function $f \in \mathbb{L}[0,1]$, define the set of points where it is nonzero:
$$
M^{\pm}_{f}:=\{t\in[0,1]:f(t)\neq 0\}.
$$

The following two theorems, which establish LDPs for $\xi_n$ and for $\sup\limits_{t\in[0,1]}|\xi_n(t)|$, are the main results of this work.

\begin{Theorem} \label{t.1}
The sequence $\xi_n$ satisfies the LDP on the metric space $(\mathbb{L}[0,1],\rho_L)$ with normalizing function $\psi(n)=n$ and rate function
\begin{equation}\label{17-05-25-1}
I(f):= \begin{cases} 
\lambda \mathfrak{m}(M^{\pm}_{f})+\frac{1}{2}\int_{M^{\pm}_{f}}\dot{f}^2(t)\,dt,  &
\text{if } f\in \widetilde{\mathbb{AC}}_0[0,1],\\
\infty, & \text{if } f\not\in \widetilde{\mathbb{AC}}_0[0,1],
\end{cases}
\end{equation}
where $\mathfrak{m}$ denotes the Lebesgue measure and $\dot{f}$ is the derivative of $f$ in the sense of absolute continuity.
\end{Theorem}

\begin{Remark}
Note that the result of Theorem~\ref{t.1} remains valid if we replace the metric space $(\mathbb{L}[0,1],\rho_L)$
by $(\mathbb{L}_p[0,1],\rho_{L_p})$, where $p\in(1,\infty)$,
$$
\mathbb{L}_p[0,1]:=\left\{f:\int_0^1|f(t)|^pdt<\infty\right\}, \quad
\rho_{L_p}(f,g):=\left(\int_0^1|f(t)-g(t)|^pdt\right)^{1/p}, \quad f,g\in\mathbb{L}_p[0,1].
$$
\end{Remark}

In addition, we establish a large deviation principle (LDP) for the maximum absolute value of the process $\xi_n(t)$ over the interval $[0,1]$, denoted as $\sup_{t \in [0,1]} |\xi_n(t)|$.

\begin{Theorem} \label{t.2} For any $x \geq 0$, the following equality holds:
\begin{equation}\label{07-06-25-2}
I_{\sup}(x):=\lim\limits_{n\rightarrow\infty}\frac{1}{n}\ln\mathbf{P}\left(\sup\limits_{t\in[0,1]}|\xi_n(t)|\geq x\right)
=\begin{cases} \infty,  &
\mbox{if }\; x<0,\\
\sqrt{2\lambda}x, & \mbox{if }\; x\in[0,\sqrt{2\lambda}],\\
\lambda+\frac{x^2}{2}, & \mbox{if }\; x>\sqrt{2\lambda}.
\end{cases}
\end{equation}
\end{Theorem}

\begin{figure}
    \centering
    \includegraphics[width=0.9\linewidth]{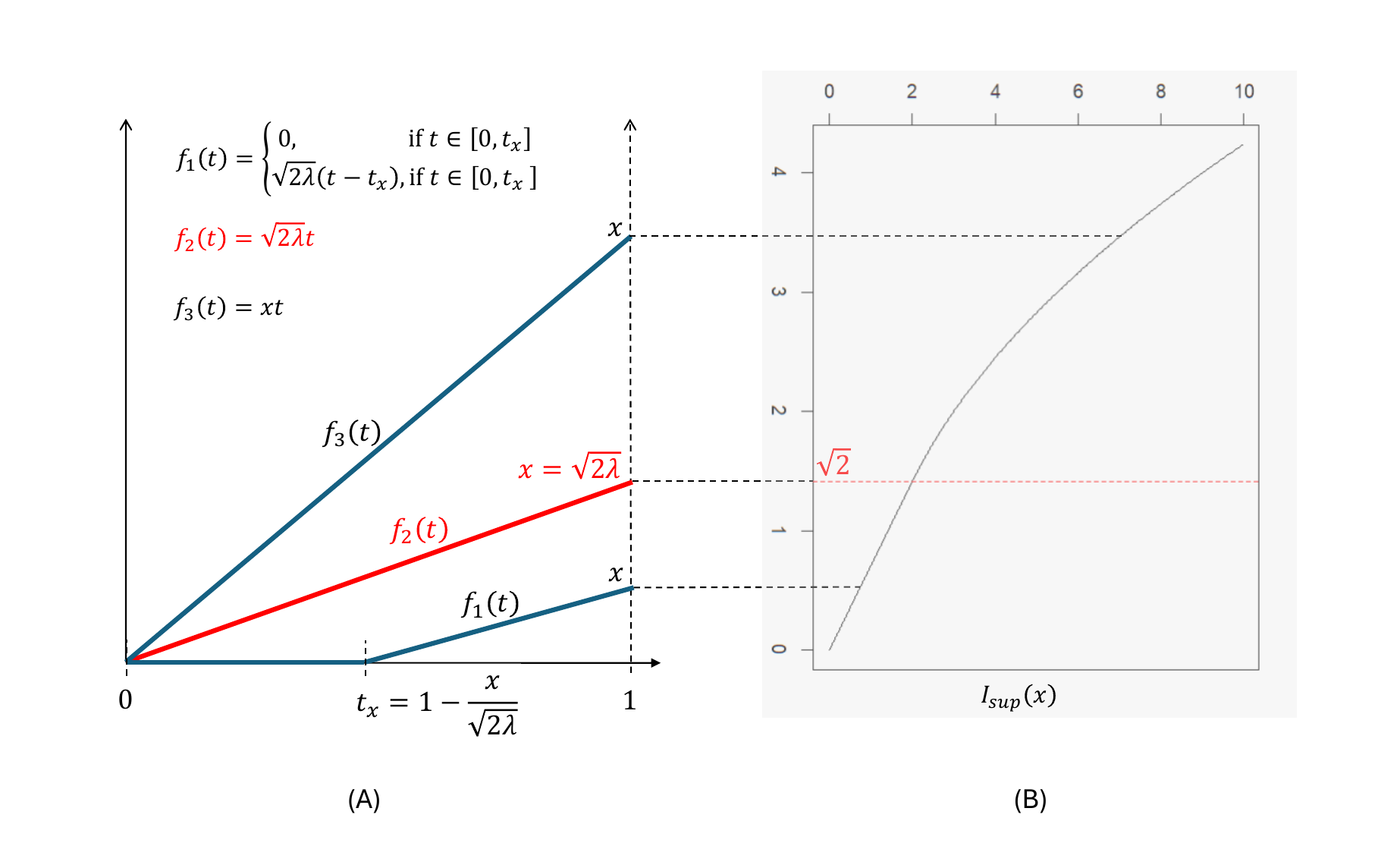}
    \caption{\textbf{Figure (A)} shows the optimal trajectories that bring the supremum of the process $\xi_n$ to the level $x$ during the time interval $[0,1]$. \textbf{Figure (B)} shows the rate function \eqref{07-06-25-2} (with values on the horizontal axis) corresponding to each level $x$ in the phase space (vertical axis), for $\lambda = 1$.}
    \label{fig3}
\end{figure}

The proof of this theorem is based on the optimization of the rate function established by Theorem~\ref{t.1}, which provides the optimal trajectories, as shown in Figure~\ref{fig3} (A). The trajectories belong to the special class of monotone trajectories \eqref{07-06-25-1}. The rate function on the phase space $\mathbb{R}_+$ is represented in Figure~\ref{fig3} (B), for the case $\lambda=1$. We have presented the graph in a slightly unconventional manner, placing the domain of the function along the vertical axis (ordinate) and the values of the rate function along the horizontal axis (abscissa). This choice allows us to connect two figures: (A) the graph of optimal trajectories reaching a given level $x$ (now represented along the horizontal axis) in the phase space, and (B) the corresponding value of the rate function. We hope that this alternative presentation will not cause discomfort to the reader and will serve as a meaningful illustration of the result obtained.

\section{Proofs of the Main Results}\label{pmr}

\subsection{Proof of Theorem~\ref{t.1}}

The proof of Theorem~\ref{t.1} follows the standard implication (see, for example, \cite[Lemma 4.1.23]{Dembo})
$$
\text{ET}\text{ and }\text{LLDP}\Rightarrow\text{LDP}.
$$
In the next section, we establish exponential tightness. As usual, the proof of the local large deviation principle (LLDP) is then divided into the upper and lower bound estimates. 

As the process $\xi_n$ can be regarded as a simple (random) transformation of the Wiener process $w_n$, see \eqref{16-05-25-1-2}, the rate function corresponding to the Wiener process naturally appears as a key element in the proofs. The following proposition characterizes the rate function, extending the result to the space $(\mathbb{L}[0,1], \rho_L)$.

\begin{Prop} \label{l.2}
The sequence of processes $w_n$ satisfies the Large Deviation Principle (LDP) in the spaces $(\mathbb{C}[0,1], \rho)$ and $(\mathbb{L}[0,1], \rho_L)$ with the normalization function $\psi(n) = n$ and the ``good" rate function $\hat{I}$
\begin{equation}\label{l.2-eq}
\hat{I}(f):= \begin{cases}\frac{1}{2}\int_0^1\dot{f}^2(t)dt,  &
\mbox{if }\; f\in \mathbb{AC}_0[0,1],\\
\infty, & \mbox{if }\; f\not\in \mathbb{AC}_0[0,1].
\end{cases}
\end{equation}
\end{Prop}
The proof of the proposition follows from \cite[Chapter 3, Theorems 2.1, 2.2]{Freidl} and auxiliary Lemma~\ref{l.1} in Section~\ref{Aux}.

\begin{Remark} \label{r.3}
We note that, in particular, the ``good" rate function of the sequence $w_n$ satisfies the Exponential Tightness (ET) property in both spaces $(\mathbb{C}[0,1], \rho)$ and $(\mathbb{L}[0,1], \rho_L)$ (see, for example, \cite[Lemma 3.5]{Feng}).
\end{Remark}

Let us start by establishing exponential tightness.

\subsubsection*{Exponential tightness}

Here we prove that the sequence $\xi_n$ is exponentially tight (ET) in the space $(\mathbb{L}[0,1],\rho_L)$.

For any positive real number $N>0$, let us define the following two sets
$$
C_N:=\{f\in \mathbb{L}[0,1]:\hat{I}(f)\leq N\}, \ \ \ C'_N:=\{f\in \mathbb{L}[0,1]:\text{Var}_{[0,1]}f\leq 2\sqrt{2N}, \ f(0)=0\},
$$
where $\text{Var}_{[0,1]}f$ denotes the total variation of the function $f$ on the interval $[0,1]$, and $\hat{I}$ is the rate function \eqref{l.2-eq}.

Note that the space $(\mathbb{L}[0,1], \rho_L)$ is complete and separable, and any probability measure on it is tight. In this setting, where a sequence of tight measures is considered, an equivalent condition for exponential tightness is known; see \cite[Lemma 3.3]{Feng}. According to this lemma, we construct a suitable expanded event to demonstrate exponential tightness by taking the $\varepsilon$-neighborhood of the set $C_N$ in the space $(\mathbb{L}[0,1], \rho_L)$.

Observe that the set $C_N$ contains only the absolutely continuous functions, $C_N\subset \mathbb{AC}_0[0,1]$, and the set $C'_N$ contains the functions with bounded variation, $C'_N\subset \mathbb{V}[0,1]$. Moreover, both sets are compact. Indeed, since $\hat{I}$ is a ``good" rate function according to Proposition~\ref{l.2}, then the set $C_N$ is a compact set. The compactness of the set $C'_N$ in $(\mathbb{L}[0,1],\rho_L)$ follows from \cite[Lemma A.1]{LY1}, which states that the set of functions $f$ with bounded variation and for some $c\geq 0$ satisfying the condition $f(0)\in[0,c]$  is compact in $(\mathbb{L}[0,1],\rho_L)$. 

For any function $f\in C_N$, using the Cauchy–Schwarz inequality, we obtain
\beq\label{28.04.5}
\text{Var}_{[0,1]} f=\int_0^1|\dot{f}(t)|dt\leq \sqrt{2 \hat{I}(f)} \le \sqrt{2N}.
\eeq


Therefore, for any $\varepsilon>0$, the inclusion holds. If a trajectory of the Wiener process $w_n$ remains uniformly within distance $\varepsilon$ of a function $f$, then the corresponding trajectory of the process $\xi_n$, obtained as a transformation \eqref{16-05-25-1-2} of $w_n$, remains uniformly within distance $2\varepsilon$ of the function $\hat{f}$, which is the image of $f$ under the same transformation \eqref{16-05-25-1-2}:
\beq\label{05.05.1}
\{\omega:\rho(w_n,f)<\varepsilon\}\subseteq\{\omega:\rho(\xi_n,\hat{f})<2\varepsilon\},
\eeq
where 
$$
\hat{f}(t):=f(t)-f\left(\frac{t_{\nu(nt)}}{n}\right), \ \ \ t\in[0,1].
$$
Here and further, $\omega$ is an element of the set $\Omega$, and the inclusions and equalities of events are understood up to events of zero probability.

It is easy to see that  $$f(0)=\hat{f}(0) \  \text{a.s.} \ \ \ \text{and} \ \ \ \text{Var}_{[0,1]}\hat{f} \leq 2\text{Var}_{[0,1]}f \ \text{a.s.}$$
Therefore, from formulas (\ref{28.04.5}), (\ref{05.05.1}) it follows that for any $\varepsilon > 0$
\beq\label{28.04.7}
\left\{\omega: w_n \in C_N^{(\varepsilon)}\right\} \subseteq \left\{\omega: \xi_n \in C_N'^{(\varepsilon)}\right\},
\eeq
where
$$
C_N^{(\varepsilon)} := \left\{f \in \mathbb{L}[0,1]: \rho(f, C_N) < \varepsilon \right\}, \ \ \
C_N'^{(\varepsilon)} := \left\{f \in \mathbb{L}[0,1]: \rho(f, C'_N) < 2\varepsilon \right\}.
$$
Since $\rho(f,g) \geq \rho_L(f,g)$, we have
\beq\label{03.05.1}
C_N'^{(\varepsilon)} \subseteq C_{N,L}'^{(\varepsilon)} := \left\{f \in \mathbb{L}[0,1]: \rho_L(f, C'_N) < 2\varepsilon \right\}.
\eeq

Combining formulas~(\ref{28.04.7}) and~(\ref{03.05.1}), the large deviation principle for $w_n$ in the metric space $(\mathbb{C}[0,1], \rho)$ (see Proposition~\ref{l.2}), the almost sure continuity of the trajectories of the process $w_n$, and the fact that the set $C_N^{(\varepsilon)} \cap \mathbb{C}[0,1]$ is open in $(\mathbb{C}[0,1], \rho)$ for any $\varepsilon > 0$ and $N > 0$, we obtain

\beq\label{30.04.11}
\begin{aligned} 
\limsup\limits_{n \rightarrow \infty} \frac{1}{n} \ln \mathbf{P} \left( \xi_n \in \bigl(C_{N,L}'^{(\varepsilon)}\bigr)^c \right)
& \leq \limsup\limits_{n \rightarrow \infty} \frac{1}{n} \ln \left(1 - \mathbf{P} \left( \xi_n \in C_N'^{(\varepsilon)} \right) \right)
\\
&\leq \limsup\limits_{n \rightarrow \infty} \frac{1}{n} \ln \left(1 - \mathbf{P} \left( w_n \in C_N^{(\varepsilon)} \right) \right) 
\\
&= \limsup\limits_{n \rightarrow \infty} \frac{1}{n} \ln \left(1 - \mathbf{P} \left( w_n \in C_N^{(\varepsilon)} \cap \mathbb{C}[0,1] \right) \right)
\\
&= \limsup\limits_{n \rightarrow \infty} \frac{1}{n} \ln \mathbf{P} \left( w_n \in \bigl(C_N^{(\varepsilon)} \cap \mathbb{C}[0,1]\bigr)^c_{(\mathbb{C}[0,1],\rho)} \right)
\\
&\leq - \inf\limits_{f \in (C_N^{(\varepsilon)} \cap \mathbb{C}[0,1])^c_{(\mathbb{C}[0,1],\rho)}} \hat{I}(f) \leq -N,
\end{aligned}
\eeq
where $\bigl(C_N^{(\varepsilon)} \cap \mathbb{C}[0,1]\bigr)^c_{(\mathbb{C}[0,1],\rho)}$
denotes the complement of the set $C_N^{(\varepsilon)} \cap \mathbb{C}[0,1]$
in the space $(\mathbb{C}[0,1], \rho)$.

The exponential tightness of the sequence $\xi_n$ follows from (\ref{30.04.11}) and \cite[Lemma 3.3]{Feng}.

\medskip

Now, let us prove the LLDP for the sequence $\xi_n$.

\subsubsection*{Upper bound}

We start by establishing the upper bound in LLDP:
\beq\label{28.04.1}
\limsup\limits_{n \rightarrow \infty} \frac{1}{n} \ln \mathbf{P}(\rho(\xi_n, f) < \varepsilon) \leq -I(f).
\eeq
Observe, by the definition of the rate function \eqref{17-05-25-1}, that $I(f)$ takes values different from $\infty$ only on functions from the set $\widetilde{\mathbb{AC}}_0[0,1]$. Recall that the definition of this set consists of three conditions; see the Introduction section.  
Indeed, if the first condition is not satisfied, then $I(f) = \infty$ because $\xi_n(0) = 0$ and the sequence $w_n$ is exponentially tight in the space $(\mathbb{C}[0,1], \rho)$ (see Remark~\ref{r.3}).  
If the second condition is not satisfied, then $I(f) = \infty$ because, by definition, the random process $\xi_n$ takes the value zero at all points of discontinuity.  
If the third condition is not satisfied, then $I(f) = \infty$ due to the large deviation principle for $w_n$ (see Proposition~\ref{l.2}) and the fact that the trajectory of $\xi_n$ between neighboring discontinuities coincides with a shifted trajectory of the process $w_n$.

Now consider $f \in \widetilde{\mathbb{AC}}_0[0,1]$.
If $f \equiv 0$, the inequality~(\ref{28.04.1}) holds trivially, since $I(f) = 0$. Therefore, we assume
that the Lebesgue measure of the set of points where $f$ does not equal zero is non-zero. For a such given $f$ and $\delta > 0$ let $s_0 := 0$ and 
$$
\begin{aligned}
  s_{2r-1} &:= 1 \wedge \inf\left\{ t \in [s_{2r-2},1]: |f(t)| = 2\delta \right\}, \\
s_{2r} &:= 1 \wedge \inf\left\{ t \in [s_{2r-1},1]: |f(t)| \leq \delta \right\}, \ \ \ r \in \mathbb{N}.  
\end{aligned}
$$
Let $m_\delta := \min\{ l : s_l = 1 \}$.
Note that $m_\delta < \infty$ since $f \in \widetilde{\mathbb{AC}}_0[0,1]$.
For definiteness, let us consider the case where $m_\delta$ is odd.
It is easy to see that by construction, the function $|f|$
is continuous and its values are no less than $\delta$ on each interval $[s_{2r-1}, s_{2r})$,
for $1 \leq r \leq \frac{m_\delta - 1}{2}$. Note also that the zeros of the function $f$ are in the intervals $[s_{2r}, s_{2r})$.

For $1 \leq r \leq \frac{m_\delta - 1}{2}$ define the events
$$
\begin{aligned} 
A_r &:= \left\{ \omega : \int_{s_{2r-1}}^{s_{2r}} |\xi_n(t) - f(t)| dt < \varepsilon \right\}, \\
B_r &:= \bigl\{ \omega : \nu(ns_{2r}-) - \nu(ns_{2r-1}) = 0 \bigr\},
\\C_r &:= \left\{ \omega : \int_{s_{2r-1}}^{s_{2r}} \left| w_n(t) - \frac{1}{n}w(t_{\nu(ns_{2r-1})}) - f(t) \right| dt \leq \varepsilon \right\}.
\end{aligned}
$$
It is easy to see that for any $\delta > 0$
on the event $A_r$, $1 \leq r \leq \frac{m_\delta - 1}{2}$, for sufficiently small $\varepsilon > 0$ and
$t \in [s_{2r-1}, s_{2r})$ the equalities 
\beq\label{28.04.8}
\begin{aligned} 
&\xi_n(t) = w_n(t) - \frac{1}{n} \sum\limits_{k=1}^{\nu(ns_{2r-1})} \xi(t_k-) = w_n(t) - \frac{1}{n} w(t_{\nu(ns_{2r-1})}), \\
&\nu(nt) - \nu(ns_{2r-1}) = 0,
\end{aligned}
\eeq
hold, i.e. on the event $A_r$, for sufficiently small $\varepsilon > 0$, the process $\xi_n(t)$ is continuous on the interval $[s_{2r-1}, s_{2r})$. In other words, the event $A_r$ implies both $C_r$ and $B_r$, separating the actions of $w_n$ and $\nu$.

Therefore, for any $\delta > 0$, $N > 0$  
and sufficiently small $\varepsilon > 0$, we have
\beq\label{28.04.2}
\begin{aligned}
\mathbf{P}(\rho_L(\xi_n,f)<\varepsilon)
& \leq\mathbf{P}\left(w_n\in K_N,~\bigcap\limits_{r=1}^{(m_\delta-1)/2}A_r\right)
+\mathbf{P}\left(w_n\in (K_N)^c\right)
\\
&\leq\mathbf{P}\left(w_n\in K_N,~\bigcap\limits_{r=1}^{(m_\delta-1)/2}C_r,~\bigcap\limits_{r=1}^{(m_\delta-1)/2}B_r\right)
+\mathbf{P}\left(w_n\in (K_N)^c\right),
\end{aligned}
\eeq
where $K_N$ is compact in metric space $(\mathbb{C}[0,1], \rho)$ such that
\beq\label{30.04.1}
\limsup\limits_{n\rightarrow\infty}\frac{1}{n}
\ln\mathbf{P}\left(w_n\in (K_N)^c \right)\leq -N.
\eeq
(The existence of such a compact set follows from Remark \ref{r.3}). Since $w_n$ and $\nu$ are independent, we aim to express the event as the intersection of two independent events, each determined by the behavior of $w_n$ (events $C_r$) and $\nu$ (events $B_r$), respectively.

Let us refine the intersection of events $C_r$ with $K_N$. Define the following closed set
$$
L_{\delta,\varepsilon}:=\bigcap\limits_{r=1}^{(m_\delta-1)/2}\left\{g\in \mathbb{L}[0,1]:\inf\limits_{-c_N\leq a_r\leq c_N}\int_{s_{2r-1}}^{s_{2r}}|g(t)-a_r-f(t)|dt\leq\varepsilon\right\},
$$
$c_N:=\sup\limits_{h\in K_N}\sup\limits_{t\in[0,1]}\,|h(t)|$.

Combining formula (\ref{28.04.2}) with the independence of the processes $w_n$ and $\nu$, as well as the independence of increments of $\nu$, we obtain that for any $\delta > 0$, $N > 0$, and sufficiently small $\varepsilon > 0$,
\beq\label{05.01.2}
\begin{aligned}
\mathbf{P}(\rho_L(\xi_n,f)<\varepsilon) &\leq
\mathbf{P}\left(w_n\in L_{\delta,\varepsilon},\bigcap\limits_{r=1}^{(m_\delta-1)/2}B_r\right)
+\mathbf{P}\left(w_n\in (K_N)^c\right)
\\
&=\mathbf{P}\left(w_n\in L_{\delta,\varepsilon}\right)\mathbf{P}\left(\bigcap\limits_{r=1}^{(m_\delta-1)/2}B_r\right)
+\mathbf{P}\left(w_n\in (K_N)^c\right)
\\
& =\mathbf{P}\left(w_n\in L_{\delta,\varepsilon}\right)\exp\left\{-\lambda n\sum\limits_{r=1}^{(m_\delta-1)/2}(s_{2r}-s_{2r-1})\right\}+\mathbf{P}\left(w_n\in (K_N)^c\right).
\end{aligned}
\eeq

It is easy to see that for any $\varepsilon > 0$, the closed set $L_{\delta,\varepsilon}$ contains a function $g$ such that $\dot{g}(t) = \dot{f}(t)$ for almost every $t \in \bigcup\limits_{r=1}^{(m_\delta - 1)/2} [s_{2r-1}, s_{2r}]$. Therefore, using Proposition~\ref{l.2} and the fact that the rate functional $\hat{I}$ is lower semicontinuous,
we obtain
\beq\label{28.04.3}
\lim\limits_{\varepsilon\rightarrow 0}\limsup\limits_{n\rightarrow\infty}\frac{1}{n}
\ln\mathbf{P}\left(w_n\in L_{\delta,\varepsilon}\right)\leq
-\lim\limits_{\varepsilon\rightarrow 0}\hat{I}(L_{\delta,\varepsilon})
\leq -\frac{1}{2}\sum\limits_{r=1}^{(m_\delta-1)/2}\int_{s_{2r-1}}^{s_{2r}}\dot{f}^2(t)dt.
\eeq
For all $\delta > 0$ and $N > 0$, inequalities (\ref{30.04.1}), (\ref{05.01.2}), and (\ref{28.04.3}) yield
$$
\begin{aligned} 
&\lim\limits_{\varepsilon\rightarrow 0}\limsup\limits_{n\rightarrow\infty}
\frac{1}{n}
\ln\mathbf{P}(\rho_L(\xi_n,f)<\varepsilon)
\\ &\leq
-\min\left(\lambda\sum\limits_{r=1}^{(m_\delta-1)/2}(s_{2r}-s_{2r-1})
+\frac{1}{2}\sum\limits_{r=1}^{(m_\delta-1)/2}\int_{s_{2r-1}}^{s_{2r}}\dot{f}^2(t)dt, N\right).
\end{aligned}
$$
Taking the limits $N \to \infty$ and then $\delta \to 0$, we obtain the inequality
$$
\begin{aligned} 
&\lim\limits_{\varepsilon\rightarrow 0}\limsup\limits_{n\rightarrow\infty}
\frac{1}{n}\ln\mathbf{P}(\rho_L(\xi_n,f)<\varepsilon)
\\
&\leq -\lambda \mathfrak{m}(\{t\in[0,1]:f(t)\neq 0\})-\frac{1}{2}\int_0^1\dot{f}^2(t)\mathbf{I}(f(t)\neq 0)dt,
\end{aligned}
$$
where $\mathbf{I}(A)$ denotes the indicator function of the set $A$. The proof of inequality (\ref{28.04.2}) for odd values of $m_\delta$ is complete. The even case is analogous and therefore omitted.

\subsubsection*{Lower bound}

Now we prove 
\beq\label{29.04.1}
\liminf\limits_{n\rightarrow\infty}\frac{1}{n}\ln \mathbf{P}(\rho_L(\xi_n,f)<\varepsilon)\geq -I(f).
\eeq

If $I(f) = \infty$, then inequality (\ref{29.04.1}) holds trivially. Therefore, we will henceforth assume that $I(f) < \infty$,  and consequently $f \in \widetilde{\mathbb{AC}}_0[0,1]$. Given such a function $f$, and following the same idea as before to construct the sequence of points, we set $s'_0 := 0$, and for $\varepsilon > 0$ we define
$$
\begin{aligned} 
s'_{2r-1}&:=1\wedge\inf\left\{t\in[s'_{2r-2},1]:|f(t)|=\frac{\varepsilon}{3}\right\}, \\
s'_{2r}&:=1\wedge\inf\left\{t\in[s'_{2r-1},1]:|f(t)|\leq\frac{\varepsilon}{4}\right\}, \ \ \ r\in \mathbb{N}.
\end{aligned}
$$
Let $m_\varepsilon := \min\{l : s'_l = 1\}$.  
Recall that $m_\varepsilon < \infty$ since $f \in \widetilde{\mathbb{AC}}_0[0,1]$.  
Consider the case when $m_\varepsilon$ is odd. 
To prove the lower bound, we define “more restrictive” events on the previously specified intervals, in the same spirit: to separate the contributions of the Poisson and Wiener processes, thereby revealing the rate function of the Wiener process.

\begin{figure}
    \centering
    \includegraphics[width=1\linewidth]{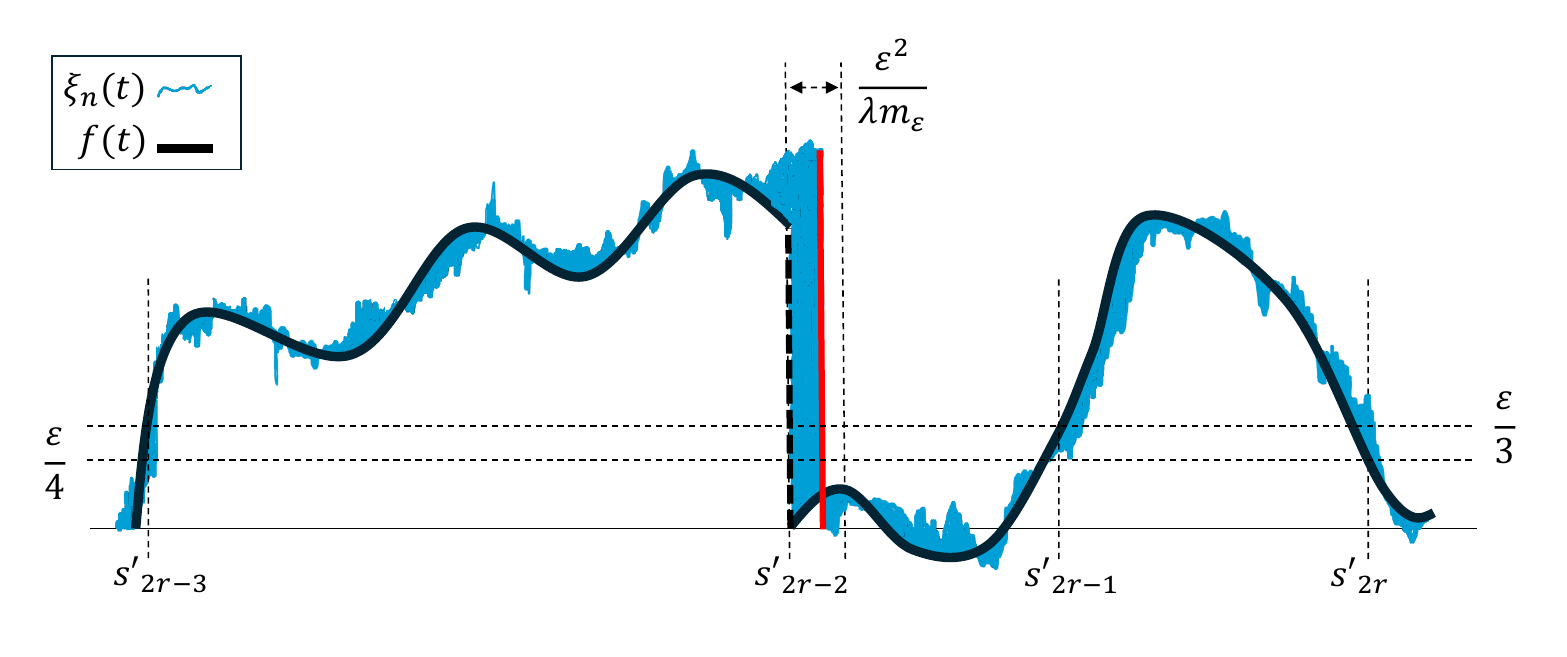}
    \caption{The figure shows the trajectory of the process $\xi_n$ (in blue) and the function $f$ (thick black line). The blue area represents the integral estimated in formulas \eqref{29.04.2} and \eqref{29.04.3}. The red line corresponds to the resetting jump of the process $\xi_n$ triggered by the event $F_r$.
}
    \label{fig2}
\end{figure}

For $1\leq r\leq \frac{m_\varepsilon+1}{2}$ we define the following events
$$
\begin{aligned}
F_r&:=\left\{\omega:\sup\limits_{t\in[s'_{2r-2},s'_{2r-1}]}\big|(\nu(nt)-\nu(ns'_{2r-2}))-\lambda n(t-s'_{2r-2})\big|<\frac{n\varepsilon^2}{m_\varepsilon}\right\}, 
\\
G_r&:=\left\{\omega:\sup\limits_{t\in[s'_{2r-2},s'_{2r-1}]}\big|w_n(t)-w_n(s'_{2r-2})\big|<\frac{\varepsilon}{8}\right\},
\end{aligned}
$$
and for $1\leq r\leq \frac{m_\varepsilon-1}{2}$ we define
$$
\begin{aligned} 
H_r&:=\left\{\omega:\sup\limits_{t\in[s'_{2r-1},s'_{2r})}\big|(w_n(t)-w_n(s'_{2r-1}))-f(t))\big|<\frac{\varepsilon}{8}\right\},
\\
Q_r&:=\Bigl\{\omega:\nu(ns'_{2r}-)-\nu(ns'_{2r-1})=0
\Bigr\}.
\end{aligned}
$$
Figure~\ref{fig2} illustrates the events defined above. Define
$$
F:=\bigcap\limits_{r=1}^{(m_\varepsilon+1)/2}F_r, \ \ \  G:=\bigcap\limits_{k=1}^{(m_\varepsilon+1)/2}G_r, \ \ \
H:=\bigcap\limits_{k=1}^{(m_\varepsilon-1)/2}H_r, \ \ \  Q:=\bigcap\limits_{k=1}^{(m_\varepsilon-1)/2}Q_r.
$$

In particular, the event $F$ ensures the occurrence of at least one jump of the process $\nu(nt)$ in the right neighborhood (of length $\frac{\varepsilon^2}{\lambda m_\varepsilon}$) of the point $s'_{2r-2}$, for every $r$ such that $1 \leq r \leq \frac{m\varepsilon + 1}{2}$.
If the event $G$ occurs, then  for all $t,s\in[s'_{2r-2},s'_{2r-1}]$, $1\leq r\leq \frac{m_\varepsilon+1}{2}$,
\beq\label{21.05.1}
\big|w_n(t)-w_n(s)\big|<\frac{\varepsilon}{4}.
\eeq  
Due to the (\ref{21.05.1}), it follows that on the event $F\cap G\cap H\cap Q$, for any $\varepsilon > 0$ and $1\leq r\leq \frac{m_\varepsilon-1}{2}$, the following inequality holds:
\beq\label{21.05.2}
\begin{aligned} 
\sup\limits_{t\in[s'_{2r-1},s'_{2r})}&\big|\xi_n(t)-f(t)\big|
=\sup\limits_{t\in[s'_{2r-1},s'_{2r})}\left|w_n(t)-\frac{1}{n}w(t_{\nu(n s'_{2r-1})})-f(t)\right|
\\
&\leq\left|w_n(s'_{2r-1})-\frac{1}{n}w(t_{\nu(n s'_{2r-1})})\right|+
\sup\limits_{t\in[s'_{2r-1},s'_{2r})}\left|w_n(t)-w_n(s'_{2r-1})-f(t)\right|
\\
&\leq\sup\limits_{t,s\in[s'_{2r-2},s'_{2r-1}]}\big|w_n(t)-w_n(s)\big|
+\sup\limits_{t\in[s'_{2r-1},s'_{2r})}\left|w_n(t)-w_n(s'_{2r-1})-f(t)\right|
\\
& <\frac{\varepsilon}{4}+\frac{\varepsilon}{8}=\frac{3\varepsilon}{8}.
\end{aligned} 
\eeq
From formulas (\ref{21.05.1}) and (\ref{21.05.2}) it follows that, on the event $F\cap G\cap H\cap Q$, for all $1\leq r\leq \frac{m_\varepsilon+1}{2}$, the following inequality holds:
\beq\label{21.05.3}
\begin{aligned} 
\sup\limits_{t,s\in[s'_{2r-2},s^*_{2r-2})}\big|\xi_n(t)\big|&<|f(s'_{2r-2}-)|
+\frac{3\varepsilon}{8}+\sup\limits_{t\in[s'_{2r-2},s^*_{2r-2})}\left|w_n(t)-w_n(s'_{2r-2})\right|
\\
&<\sup\limits_{t\in[0,1]}|f(t)|+\frac{3\varepsilon}{8}+\frac{\varepsilon}{4} < \sup\limits_{t\in[0,1]}|f(t)|+1,
\end{aligned}
\eeq
where $s^*_{2r-2}$ denotes the first jump time of the process $\nu(nt)$ within the interval $[s'_{2r-2},s'_{2r-1})$.

Applying formulas (\ref{21.05.1}) and (\ref{21.05.3}), we obtain that on the event $F\cap G\cap H\cap Q$, for all sufficiently small $\varepsilon>0$,
\beq\label{29.04.2}
\begin{aligned} 
\sum\limits_{r=1}^{(m_\varepsilon+1)/2}\int_{s'_{2r-2}}^{s'_{2r-1}}&|\xi_n(t)-f(t)|dt \\
& < 
\left(2\sup\limits_{t\in[0,1]}|f(t)|+1\right)
\cdot\frac{\varepsilon^2}{\lambda m_\varepsilon}\cdot\frac{m_\varepsilon+1}{2}+
\frac{7\varepsilon}{12}\sum\limits_{r=1}^{(m_\varepsilon+1)/2}(s'_{2r-1}-s'_{2r-2}) 
\\
& < \frac{2\sup\limits_{t\in[0,1]}|f(t)|+1}{\lambda}\varepsilon^2
+\frac{7\varepsilon}{12}\sum\limits_{r=1}^{(m_\varepsilon+1)/2}(s'_{2r-1}-s'_{2r-2})< \frac{\varepsilon}{24} + \frac{7\varepsilon}{12} = \frac{5\varepsilon}{8}.
\end{aligned}
\eeq
It follows from formula (\ref{21.05.2}) that, on the event $F\cap G\cap H\cap Q$,

\beq\label{29.04.3}
\sum\limits_{r=1}^{(m_\varepsilon-1)/2}\int_{s'_{2r-1}}^{s'_{2r}}|\xi_n(t)-f(t)|dt
<\sum\limits_{r=1}^{(m_\varepsilon-1)/2}\int_{s'_{2r-1}}^{s'_{2r}}\frac{3\varepsilon}{8}dt
\leq\frac{3\varepsilon}{8},
\eeq 
From formulas (\ref{29.04.2}) and (\ref{29.04.3}) it follows that for sufficiently small $\varepsilon > 0$,  
on the event $F \cap G \cap H \cap Q$, the inequality $\rho_L(\xi_n, f) < \varepsilon$ holds.  
Therefore, using the independence of the processes $\nu$ and $w_n$, as well as the independence of their increments, we obtain  
for any sufficiently small $\varepsilon > 0$:
\beq\label{29.04.5}
\begin{aligned} 
\mathbf{P}(\rho_L(\xi_n,f)<\varepsilon)
&\geq\mathbf{P}\left(F,G,H,Q\right) =
\mathbf{P}\left(F\right)
\mathbf{P}\left(G\right)
\mathbf{P}\left(H\right)
\mathbf{P}\left(Q\right)
\\
&
=\mathbf{P}\left(F\right)
\mathbf{P}\left(G\right)
\mathbf{P}\left(H\right)
\exp\left\{-\lambda n\sum\limits_{r=1}^{(m_\varepsilon-1)/2}(s'_{2r}-s'_{2r-1})\right\}.
\end{aligned}
\eeq
Denote $v:=\lambda t$, $t\in[0,1]$. Applying Lemma~\ref{l.3}, we obtain fro any $\varepsilon>0$
\beq\label{29.04.7}
\liminf\limits_{n\rightarrow\infty}\frac{1}{n}\ln\mathbf{P}\left(F\right)
\geq\liminf\limits_{n\rightarrow\infty}\frac{1}{n}\ln\mathbf{P}\left(\sup\limits_{t\in[0,1]}|\nu_n(t)-\lambda t|
<\frac{\varepsilon^2}{2m_\varepsilon}\right)\geq -I_1(v)=0.
\eeq
According to Proposition~\ref{l.2}, for any $\varepsilon>0$ we have
\beq\label{29.04.8}
\liminf\limits_{n\rightarrow\infty}\frac{1}{n}\ln\mathbf{P}\left(G\right)
\geq\liminf\limits_{n\rightarrow\infty}\frac{1}{n}\ln\mathbf{P}\left(\rho(w_n,0)
<\frac{\varepsilon}{8}\right)\geq -\hat{I}(0)=0.
\eeq
Consider the following open set of functions
$$
C_f:=
\left\{g\in \mathbb{C}[0,1]:\sup\limits_{t\in[s'_{2r-1},s'_{2r})}
\big|(g(t)-g(s'_{2r-1}))-f(t)\big|<\frac{\varepsilon}{8}, \ 1\leq r\leq \frac{m_\varepsilon-1}{2}
\right\}.
$$
Proposition~\ref{l.2} yields that for any $\varepsilon>0$
\beq\label{29.04.10}
\liminf\limits_{n\rightarrow\infty}\frac{1}{n}\ln\mathbf{P}\left(H\right)
=\liminf\limits_{n\rightarrow\infty}\frac{1}{n}\ln\mathbf{P}\left(w_n\in C_f\right)\geq -\hat{I}(C_f)\geq
-\frac{1}{2}\sum\limits_{r=1}^{(m_\varepsilon-1)/2}\int_{s'_{2r-1}}^{s'_{2r}}\dot{f}^2(t)dt.
\eeq
Utilizing (\ref{29.04.5})--(\ref{29.04.10}), for a sufficiently small $\varepsilon>0$ we obtain
$$
\liminf\limits_{n\rightarrow\infty}\frac{1}{n}\ln\mathbf{P}(\rho_L(\xi_n,f)<\varepsilon)
\geq-\lambda\sum\limits_{k=1}^{(m_\varepsilon-1)/2}(s'_{2k}-s'_{2k-1})
-\frac{1}{2}\sum\limits_{k=1}^{(m_\varepsilon-1)/2}\int_{s'_{2k-1}}^{s'_{2k}}\dot{f}^2(t)dt.
$$
Taking the limit as $\varepsilon \rightarrow 0$ completes the proof of the formula (\ref{29.04.1}) when $m_\varepsilon$ is odd. The proof for even $m_\varepsilon$ is analogous, so we omit it.

Formulas (\ref{28.04.1}) and (\ref{29.04.1}) establish the LLDP for the sequence $\xi_n$, and exponential tightness completes the proof of Theorem~\ref{t.1}. $\Box$


\subsection{Proof of Theorem~\ref{t.2}}

The proof is divided into the upper and lower bounds. We begin with the proof of the upper bound
\beq\label{02.05.3}
\limsup\limits_{n\rightarrow\infty}\frac{1}{n}\ln\mathbf{P}\left(\sup\limits_{t\in[0,1]}|\xi_n(t)|\geq x\right)
\leq-I_{\sup}(x).
\eeq
Define the following set of functions
$$
A_x:=\left\{f\in\mathbb{L}[0,1]:\underset{t\in[0,1]}{\text{ess sup}}\,|f(t)|\geq x\right\}.
$$
Since the process $\xi_n$ has sample paths in the space $\mathbb{D}[0,1]$ almost surely, its supremum coincides almost surely with its essential supremum, i.e.,
$$
\left\{\omega:\sup\limits_{t\in[0,1]}|\xi_n(t)|\geq x\right\}
=\left\{\omega:\xi_n\in A_x\right\}.
$$
Note that the upper bound in the LDP from Theorem~\ref{t.1} cannot be applied directly to the set $A_x$, since it is not closed in the metric space $(\mathbb{L}[0,1],\rho_L)$. However, any set consisting of functions with a maximum of at least $x$ and uniformly continuous at their maximum points is closed. We denote this set by $B_x$. By the Ascoli–Arzel\`a theorem and the definition of the process $\xi_n$, for any compact set $K'_N$ in the space $(\mathbb{C}[0,1],\rho)$ we have
$$
\left\{\omega:\xi_n\in B_x\right\}\supseteq \left\{\omega:\xi_n\in A_x,w_n\in K'_N\right\}.
$$
We assume that the compact set $K'_N$ is chosen in such a way that inequality (\ref{30.04.1}) holds for the probability of its complement. Using the notation introduced above, we obtain
\beq\label{02.05.1}
\begin{aligned} 
\mathbf{P}\left(\sup\limits_{t\in[0,1]}|\xi_n(t)|\geq x\right)
&\leq\mathbf{P}\left(\xi_n\in A_x,w_n\in K'_N\right)+
\mathbf{P}\left(w_n\in (K'_N)^c\right)
\\
&\leq\mathbf{P}\left(\xi_n\in B_x\right)+\mathbf{P}\left(w_n\in (K'_N)^c\right).
\end{aligned}
\eeq
Inequality (\ref{02.05.1}) together with Theorem~\ref{t.1} implies that for any $N>0$
$$
\limsup\limits_{n\rightarrow\infty}\frac{1}{n}\ln\mathbf{P}\left(\sup\limits_{t\in[0,1]}|\xi_n(t)|\geq x\right)
\leq-\min(I(B_x),N)\leq -\min(I(A_x),N).
$$
Taking the limit $N \rightarrow \infty$, we obtain
$$
\limsup\limits_{n\rightarrow\infty}\frac{1}{n}\ln\mathbf{P}\left(\sup\limits_{t\in[0,1]}|\xi_n(t)|\geq x\right)
\leq -\inf\limits_{f\in A_x}I(f).
$$
From the definition of the functional $I$, it follows that the infimum is attained over functions in the set $A_x \cap \widetilde{\mathbb{AC}}_0[0,1]$.

It is easy to see that for $f \in A_x \cap \widetilde{\mathbb{AC}}_0[0,1]$, the following inequality holds:
$$
x\leq\sup\limits_{t\in[0,1]}|f(t)|\leq\int_{M_f^{\pm}}|\dot{f}(t)|dt.
$$
Therefore, by the Cauchy–Schwarz inequality, we obtain
\beq\label{02.05.2}
x^2\leq \mathfrak{m}(M^{\pm}_{f})\int_{M_f^{\pm}}\dot{f}^2(t)dt.
\eeq
The inequality (\ref{02.05.2}) yields
$$
\begin{aligned} 
 \inf\limits_{f\in A_x\cap\widetilde{\mathbb{AC}}_0[0,1]}
 I(f)
&\geq\inf\limits_{f\in A_x\cap\widetilde{\mathbb{AC}}_0[0,1]}
\left(\lambda\mathfrak{m}(M^{\pm}_{f})+\frac{1}{2}\int_{M_f^{\pm}}\dot{f}^2(t)dt\right)
\\ &\geq
\inf\limits_{\mathfrak{m}(M^{\pm}_{f})\in [0,1]}
\left(\lambda\mathfrak{m}(M^{\pm}_{f})+ \frac{x^2}{2\mathfrak{m}(M^{\pm}_{f})}\right).
\end{aligned}
$$
Applying standard methods of differential calculus completes the proof of formulas (\ref{02.05.3}).

\vspace{0.5cm}
Now we prove the lower bound
\beq\label{02.05.5}
\liminf\limits_{n\rightarrow\infty}\frac{1}{n}\ln\mathbf{P}\left(\sup\limits_{t\in[0,1]}|\xi_n(t)|\geq x\right)
\geq-I_{\sup}(x).
\eeq

Using Theorem~\ref{t.1} and the fact that the set
$$
C_x:=\left\{f\in\mathbb{L}[0,1]:\underset{t\in[0,1]}{\text{ess sup}}\,|f(t)|> x\right\}
$$
is the open set in the metric space $(\mathbb{L}[0,1],\rho_L)$, we obtain
\beq\label{02.05.7}
\liminf\limits_{n\rightarrow\infty}\frac{1}{n}\ln\mathbf{P}\left(\sup\limits_{t\in[0,1]}|\xi_n(t)|\geq x\right)
\geq
\liminf\limits_{n\rightarrow\infty}\frac{1}{n}\ln\mathbf{P}\left(\xi_n\in C_x\right)
\geq-\inf\limits_{f\in C_x}I(x).
\eeq
It is easy to see that
\begin{equation}\label{07-06-25-1}
\left\{f\in\mathbb{L}[0,1]:f(t)=0,~t\in[0,s];~f(t)=k(t-s),~t\in(s,1], \text{ where } k>\frac{x}{1-s}\right\}
\subset C_x.
\end{equation}
Thus,
\beq\label{02.05.8}
\inf\limits_{f\in C_x}I(x)\leq \inf\limits_{s\in [0,1],~k>\frac{x}{1-s}}
\left(\lambda\mathfrak{m}([s,1])+\frac{k^2}{2}(1-s)\right)
=\inf\limits_{s\in [0,1]}
\left(\lambda(1-s)+\frac{x^2}{2(1-s)}\right).
\eeq
Using equality (\ref{02.05.8}) and standard methods of differential calculus, we obtain the equality
\beq\label{02.05.10}
\inf\limits_{f\in C_x}I(x)\leq\inf\limits_{s\in [0,1]}
\left(\lambda(1-s)+\frac{x^2}{2(1-s)}\right)=I_{\sup}(x).
\eeq
Formula (\ref{02.05.5}) follows from equalities (\ref{02.05.7}) and (\ref{02.05.10}). It is easy to see that the proof of Theorem \ref{t.2} mainly relies on equalities (\ref{02.05.3}) and (\ref{02.05.5}).
$\Box$

\section{Auxiliary Results}\label{Aux}

In this section, we prove some auxiliary lemmas we used to prove the theorems.

\begin{Lemma} \label{l.1}
The sequence of random processes $\eta_n$ satisfies an LDP in the space $(\mathbb{C}[0,1],\rho)$ with speed $\psi(n)$ and the "good" rate function $\tilde{I}(f)$.  
Then the sequence $\eta_n$ also satisfies an LDP in the space $(\mathbb{L}[0,1],\rho_L)$ with the same speed and rate function.
\end{Lemma}

\noindent
\DD. Consider the operator 
$
\mathcal{A} : (\mathbb{C}[0,1], \rho) \to (\mathbb{L}[0,1], \rho_L)$,
defined by the identity restriction $\mathcal{A}f = f$. Since 
$\rho(f,g) \geq \rho_L(f,g), \quad f,g \in \mathbb{C}[0,1],$
the operator $\mathcal{A}$ is continuous. Therefore, applying the contraction principle (see, e.g., \cite[Theorem 4.2.1]{Dembo}), we obtain
$$
I_L(f) = \inf_{g : \mathcal{A}g = f} \tilde{I}(g) = \tilde{I}(f),
$$
where $I_L$ is the rate function for the LDP of the sequence $\eta_n$ in the space $(\mathbb{L}[0,1], \rho_L)$.
$\Box$

\vspace{0.5cm}

Let $\nu_n(t):=\frac{\nu(nt)}{n}$, $t\in[0,1]$.

\begin{Lemma} \label{l.3} The sequence of random processes $\nu_n$ satisfies the LDP in the metric space $(\mathbb{D}[0,1],\rho)$  with normalized function $\psi(n)=n$ and ``good"  rate function
$$
I_1(f):= \begin{cases}
\int_0^1\left(\dot{f}(t)\ln\left(\frac{\dot{f}(t)}{\lambda}\right)-\dot{f}(t)+\lambda\right)dt,  &
\mbox{if }\; f\in \mathbb{AC}_{0}^{\uparrow}[0,1],\\
\infty, & \mbox{if }\; f\not\in \mathbb{AC}_{0}^{\uparrow}[0,1],
\end{cases}
$$
Here, $\mathbb{AC}_{0}^{\uparrow}[0,1]$ denotes the set of monotone non-decreasing absolutely continuous functions defined on the interval $[0,1]$ such that $f(0) = 0$.
\end{Lemma}

\noindent
\DD. This follows from \cite{Mog1} and, for example, Lemma 4.4 in \cite{LY2}. $\Box$

\vspace{0.5cm}
Note that the LDP for the stochastic process $\nu_n$ in the space $(\mathbb{D}[0,1],\rho)$ is formulated for the sequence of probability measures $\mathbf{P}(\nu_n \in A)$, not on the Borel $\sigma$-algebra $\mathfrak{B}(\mathbb{D}[0,1])$ generated by all open subsets of the metric space $(\mathbb{D}[0,1],\rho)$, but rather on the $\sigma$-algebra generated by cylinder sets defined through finite-dimensional distributions. This stems from the fact that not all sets $A \in \mathfrak{B}(\mathbb{D}[0,1])$ are measurable (for more details, see \cite[Chapter 3, Section 15]{Bil}).

\end{document}